\documentclass{amsart}

\usepackage{enumerate}

\newtheorem{theorem}{Theorem}
\newtheorem{proposition}[theorem]{Proposition}
\newtheorem{lemma}[theorem]{Lemma}

\begin{document}

\title[On minimal Rolle's domains]
      {On minimal Rolle's domains \\ for complex polynomials}
\author{Michael J. Miller}
\address{Department of Mathematics\\Le Moyne College\\Syracuse, New 
York 13214\\USA}
\email{millermj@lemoyne.edu}

\newcommand{\polys}{\mathcal P}

\begin{abstract}
Define a subset of the complex plane to be a \emph{Rolle's domain} if it
contains (at least) one critical point of every complex polynomial P
such that $P(-1)=P(1)$.  Define a Rolle's domain to be \emph{minimal}
if no proper subset is a Rolle's domain.  In this paper, we investigate
minimal Rolle's domains.
\end{abstract}

\maketitle

For real-valued functions defined on the reals, the classical Rolle's
theorem implies that if $f$ is differentiable everywhere and
$f(-1)=f(1)$ then $f$ has a critical point in the interval $(-1,1)$.
This statement is not true for general complex functions, as is
demonstrated by $f(z)=e^{i \pi z}$.

In this paper, we extend Rolle's theorem by restricting our attention
to complex polynomials.  Define $\polys$ to be the set of all complex
polynomials $P$ with $P(-1)=P(1)$.  In \cite[Section 1.1]{Sendov},
Sendov defines a \emph{Rolle's domain} to be a subset of the complex
plane that contains (at least) one critical point of every polynomial
in $\polys$.  In this paper, we define such a Rolle's domain $R$ to be
\emph{minimal} if no proper subset of~$R$ is a Rolle's domain.

An example of a Rolle's domain that is almost minimal is given by
\begin{theorem}\label{theorem1} 
The right half-plane $H=\{z: \Re(z) \ge 0\}$ is a Rolle's domain, but
is not minimal.  However, any subset of $H$ that is a Rolle's domain
must contain $\{z: \Re(z) > 0\} \cup I$, where $I=\{i \cdot \cot(k
\pi/n) : n \ge 2 \text{ and } 1 \le k< n \}$ is a dense subset of the
imaginary axis.
\end{theorem}

To prove this theorem, we will need
\begin{lemma}\label{lemma2} 
The roots $\zeta$ of the polynomial $\int_{-1}^1 (w-\zeta)^{n-1} dw$ are given by
$\zeta=i \cot(k\pi/n)$ for $1 \le k < n$.
\end{lemma}
\begin{proof}
The polynomial $\int_{-1}^1 (w-\zeta)^{n-1} dw =
\frac{(-1)^n}{n}[(\zeta-1)^n-(\zeta+1)^n]=0$ when
$(\zeta-1)/(\zeta+1)=e^{2\pi i k/n}$ for some $k$ with $1 \le k < n$.
Letting $w=e^{\pi i k/n}$, then $|w|=1$ implies that 
\begin{displaymath}
\zeta = \frac{1+w^2}{1-w^2} = \frac{2 i \Im(w^2)}{2-2 Re(w^2)} 
      = i \frac{\Re(w)}{\Im(w)}
\end{displaymath}
and the conclusion follows.
\end{proof}

We now begin the

\begin{proof}[Proof of Theorem \ref{theorem1}]
The Grace-Heawood Theorem \cite[Supplement to
  Theorem~4.3.1]{Rahman-Schmeisser} states that if $P \in \polys$ is of
degree $n$, then every circular domain containing all the points
$I_n=\{i\cot(k\pi/n): 1\le k <n\}$ must contain at least one critical
point of $P$.  Since $I_n \subset H$, then $H$ is a Rolle's domain.

Choose any irrational $\alpha \in (0,1)$ and let
$\zeta=i\cot(\alpha\pi)$.  Take any $P \in \polys$ of degree~$n$, and
note that $\zeta\notin I_n$.  Construct the circles
\begin{displaymath}
\{z : |z-(\zeta-m)|=m+\epsilon_m\} \text{ for $m=1, 2, \dots$}
\end{displaymath}
where each $\epsilon_m>0$ is chosen small enough that the points of
$I_n$ are all outside the circle.  By the Grace-Heawood theorem, the
exterior of each of these circles must contain a critical point of
$P$.  Since every point not in $H-\{\zeta\}$ is eventually inside
these circles for all sufficiently large $m$, then $P$ must have a
critical point in $H-\{\zeta\}$.  Thus $H-\{\zeta\}$ is a Rolle's
domain, and so $H$ is not minimal.

Now let $S$ be any subset of $H$ that is a Rolle's domain, and take any
$\zeta \in \{z: \Re(z) > 0\} \cup I$. If $\Re(\zeta)>0$ then
$\Re(-1/(3\zeta))<0$ so $-1/(3\zeta) \notin S$.  Since $P(z)=2\zeta
z^3+(1-3\zeta^2)z^2-2\zeta z \in \polys$ has critical points at
$\zeta$ and $-1/(3\zeta)$, then $\zeta \in S$.  If instead $\zeta \in I$ then
$\zeta=i \cdot \cot(k \pi/n)$ for some $n \ge 2$ and $1 \le k< n$.
Define $P(z)=\int_{-1}^z (w-\zeta)^{n-1} dw$ and note that by Lemma
\ref{lemma2}, we have $P \in \polys$.  Since $P$ has a critical point
only at $\zeta$, then $\zeta \in S$.  Thus we know that that $S$ must
contain $\{z: \Re(z) > 0\} \cup I$.

Finally, note that $\{k/n : n \ge 2 \text{ and } 1 \le k< n \}$ is a
dense subset of $(0,1)$ and that $t \to i \cot(t\pi)$ maps $(0,1)$
onto the imaginary axis, so $I$ is a dense subset of the imaginary
axis.
\end{proof}

In \cite[Statements 4 and 5]{Sendov}, Sendov shows that any Rolle's
domain that is symmetric with respect to both the real and imaginary
axes must contain both the sets $\{z:|\Im(z)|>1/\pi\}$ and
$\{z:|z|<1\}$, and conjectures  \cite[Conjecture 2]{Sendov} that the union
\begin{displaymath}
\{z:|\Im(z)|>1/\pi\} \cup \{z:|z|<1\}
\end{displaymath}
is a Rolle's domain (and thus the smallest symmetric Rolle's
domain).  We provide a counterexample to this conjecture with
\begin{theorem}\label{theorem3} If $c \ge 0$ and the set
$\{z:|\Im(z)|>c\} \cup \{z:|z|<1\}$
is a Rolle's domain, then $c=0$.
\end{theorem}

To prove this theorem we will need
\begin{lemma}\label{lemma4}\ 
If $0 \le i \le 5$ then
\begin{equation*}
\lim_{n \to \infty} n^{6-i} \sum_{k=i}^5 
       \binom{5}{k} \binom{k}{i} \frac{(-1)^k}{n+6-k} = \frac{-120}{i!}.
\end{equation*}
\end{lemma}
\begin{proof}
Note that $\binom{m+1}{j}= \binom{m}{j}+
\binom{m}{j-1}$, and from this a straightforward induction on $m$
verifies that
\begin{equation*}
\sum_{j=0}^m \binom{m}{j} \frac{(-1)^j}{n+1+j} 
         = \frac{m!}{\prod_{k=1}^{m+1} (n+k)} \text{\quad for $m\ge0$.}
\end{equation*}
Then since $ \binom{5}{k} \binom{k}{i} = \binom{5}{i} \binom{5-i}{5-k}$,
we have
\begin{equation*}
\begin{split}
  n^{6-i} \sum_{k=i}^5 \binom{5}{k} \binom{k}{i} \frac{(-1)^k}{n+6-k} 
 &= n^{6-i} \binom{5}{i} \sum_{k=i}^5 \binom{5-i}{5-k} \frac{(-1)^k}{n+6-k} \\
 &= -n^{6-i} \binom{5}{i} \sum_{j=0}^{5-i} \binom{5-i}{j} \frac{(-1)^j}{n+1+j}\\
 &= -n^{6-i} \binom{5}{i} \frac{(5-i)!}{\prod_{k=1}^{6-i} (n+k)}\\
 &= -\frac{5!}{i!} \frac{n^{6-i}}{\prod_{k=1}^{6-i} (n+k)}
\end{split}
\end{equation*}
and the result follows.
\end{proof}

We will also need
\begin{proposition}\label{prop5}
Given any $\epsilon>0$, for each sufficiently large $n$ there exists a
complex number $\zeta$ such that $\Re(\zeta)>1$ and
$|\Im(\zeta)|<\epsilon$ and
\begin{equation}\label{eq1}
\int_{-1}^1  (w+1)^n (w-\zeta)^5 \,dw = 0.
\end{equation}
\end{proposition}
\begin{proof}
By the binomial theorem we have
$(w-\zeta)^5 = \sum_{k=0}^5 (-1)^k \binom{5}{k}(w+1)^{5-k}(\zeta+1)^k$
so for every positive integer $n$ we have
\begin{equation*}
\begin{split}
\int_{-1}^1 (w+1)^n(w-\zeta)^5\,dw
   &= \sum_{k=0}^5 (-1)^k \binom{5}{k}(\zeta+1)^k 
               \int_{-1}^1 (w+1)^{n+5-k} \,dw\\
   &= \sum_{k=0}^5 (-1)^k \binom{5}{k}(\zeta+1)^k 
               \frac{ 2^{n+6-k}}{n+6-k}.
\end{split}
\end{equation*}
Define
\begin{equation*}
   f(z) = \sum_{k=0}^5 \binom{5}{k} 
               \frac{ 2^{n+6-k}}{n+6-k}\ z^k \text{\quad and \quad} 
   g(w) = \frac{n^6}{2^{n+6}} 
               f\left(2\left(\frac{w}{n}-1\right)\right).
\end{equation*}
Then
\begin{equation*}
\begin{split}
   g(w) &= \frac{n^6}{2^{n+6}} \sum_{k=0}^5 \binom{5}{k}
           \frac{ 2^{n+6-k}}{n+6-k} 2^k \left(\frac{w}{n}-1\right)^k \\
        &= n^6 \sum_{k=0}^5 \binom{5}{k} \frac{1}{n+6-k}
           \sum_{i=0}^k \binom{k}{i}  \left(\frac{w}{n}\right)^i (-1)^{k-i}\\
        &= \sum_{i=0}^5 (-1)^i w^i 
          \left[ n^{6-i} \sum_{k=i}^5 \binom{5}{k} \binom{k}{i} 
                                                  \frac{(-1)^k}{n+6-k}\right].
\end{split}
\end{equation*}
Then by Lemma \ref{lemma4} we know that $\lim_{n \to \infty} g(w) =
120 \sum_{i=0}^5 (-1)^{i+1} w^i/i!$, which has a root at approximately
$-0.24+3.13i$.  By \cite[Theorem 1.3.1]{Rahman-Schmeisser} the roots
of a polynomial are continuous functions of its coefficients, so for
each sufficiently large $n$ there is a root $w_0$ of $g(w)$ with
$\Re(w_0)<0$ and $0<\Im(w_0)<4$.  Require in addition that
$n>8/\epsilon$, and define $\zeta = 1-2w_0/n$.  Then
$-(\zeta+1)=2(\frac{w_0}{n}-1)$ is a root of $f$, so $\zeta$ is a root
of equation \eqref{eq1} and
\begin{equation*}
   \Re(\zeta) = 1-\frac{2}{n}\Re(w_0)>1  \text{\quad and \quad}
   |\Im(\zeta)| = \frac{2}{n}|\Im(w_0)| < 8/n < \epsilon.
\end{equation*}
\end{proof}

We can now write the
\begin{proof}[Proof of Theorem \ref{theorem3}]
Choose any $\epsilon>0$, and for a sufficiently large value of $n$
choose~$\zeta$ as in Proposition \ref{prop5}.  Define
the polynomial
\begin{displaymath}
P(z)=\int_{-1}^z (w+1)^n (w-\zeta)^5\,dw.
\end{displaymath}
Then $P(1)=0=P(-1)$ so $P \in \polys$, and $P$ has critical points
only at $-1$ and $\zeta$, so every Rolle's domain must contain either
$-1$ or $\zeta$.  Note that the set $\{z: |z|<1\}$ contains neither,
and the set $\{z: |\Im(z)|>c\}$ does not contain $-1$, so it must
contain $\zeta$.  Thus $0 \le c<\epsilon$ for every $\epsilon>0$, and
the result follows.

\end{proof}

Theorem \ref{theorem3} states that for $c \ge 0$ the only potential
Rolle's domain of the form
\begin{equation*}
\{z:|\Im(z)|>c\} \cup \{z:|z|<1\}
\end{equation*}
has $c=0$.  This region is a Rolle's domain when $c=0$, as referenced
by Sendov in \cite[Theorem 4]{Sendov} and proved in \cite[Theorem
  4.3.4]{Rahman-Schmeisser}. We present a more elementary proof in
\begin{theorem}\label{theorem6} The set $R=\{z:\Im(z)\ne0\} \cup 
\{z:|z|<1\}$ is a Rolle's domain.
\end{theorem}
\begin{proof}
Take any $P \in \polys$.  We will prove by contradiction that $P$ has
a critical point in $R$.

Assume (without loss of generality) that $P$ is monic.  If $P$ has no
critical points in $R$, then all critical points of $P$ are in the
union of intervals $(-\infty, -1] \cup [1, \infty)$, so $P'$ is a real
polynomial.  Define $Q(z)=\int_0^z P'(w)\,dw$ and note that $Q$ is
a real polynomial and that $Q(z)=P(z)-P(0)$.  Then $Q(1)=Q(-1)$ so
by Rolle's theorem for real polynomials we know that $Q$ (and thus
$P$) has a critical point in the interval $(-1,1)$, which is a
contradiction.
\end{proof}

The author thanks William Calbeck for pointing out a flaw in a
previous version of Theorem \ref{theorem1}.



\begin{thebibliography}{99}

\bibitem{Rahman-Schmeisser}
Q. I. Rahman and G. Schmeisser,
\emph{Analytic theory of polynomials},
London Mathematical Society Monographs (New Series), \textbf{26}, 
Oxford University Press, Oxford, 2002,
MR~1954841 (2004b:30015).

\bibitem{Sendov}
Bl. Sendov,
\emph{Complex analogues of the Rolle's theorem},
Serdica Math. J. \textbf{33} (2007), 387--398,
MR~2418186 (2009c:30020).

\end{thebibliography}
\end{document}